
\magnification=\magstep1
\baselineskip=14pt


\font\tenmsb=msbm10
\font\sevenmsb=msbm7
\font\fivemsb=msbm5
\newfam\msbfam
\textfont\msbfam=\tenmsb
\scriptfont\msbfam=\sevenmsb
\scriptscriptfont\msbfam=\fivemsb
\def\Bbb#1{{\fam\msbfam\relax#1}}
\font\teneufm=eufm10
\font\seveneufm=eufm7
\font\fiveeufm=eufm5
\newfam\eufmfam
\textfont\eufmfam=\teneufm
\scriptfont\eufmfam=\seveneufm
\scriptscriptfont\eufmfam=\fiveeufm

\def\C{{\Bbb C}}

\def\N{{\Bbb N}}
\def\P{{\Bbb P}}
\def\Q{{\Bbb Q}}
\def\R{{\Bbb R}}

\def\Z{{\Bbb Z}}


\def\nnproclaim#1
  {\medbreak
   \smallskip
   \noindent
   {\bf#1}\
   \begingroup\sl\penalty80}

\def\endproclaim{\endgroup\smallskip\goodbreak}

\def\proof{
  \noindent
  {\bf Proof:}
}


\def\res{\mathop{{\rm Res}}\nolimits}

\def\ini{\mathop{{\rm in}}\nolimits}

\def\ip#1#2{\langle #1, #2\rangle}

\magnification=1200
\parindent=0pt
\parskip=6pt
\font\bigbf=cmbx10 scaled \magstep1 
\def\gb{Gr\"obner basis }
\def\x{{\bf x}}
\def\y{{\bf y}}
\def\g{{\bf g}}
\def\ii{{\bf i}}
\def\jj{{\bf j}}
\def\w{{\bf w}}
\def\d{{\bf d}}
\def\y{{\bf y}}
\def\a{{\bf a}}
\def\b{{\bf b}}
\def\c{{\bf c}}
\def\r{{\bf r}}
\def\rr{{\bf r+1}}
\def\k{{\bf k}}

\def\K{{\bf K}}
\def\OO{{\cal O}}
\def\II{{\cal I}}
\def\NF{{\cal NF}}

\def\ip#1#2{\langle #1, #2\rangle}
\def\dw#1{\deg_{\w}(#1)}
\def\wp{\P^n_{\w}}
\def\iw#1{{\rm in}_{\w}(#1)}
\def\ges{g_1\ldots g_n}
\def\fes{f_1\ldots f_n}
\def\cges{g_1,\ldots, g_n}
\def\cfes{f_1,\ldots, f_n}
\def\rg#1{{\rm Res}_{{\bf g}}(#1)}
\def\rh{\displaystyle{{\rm Res}\left(
{{h\, d\x}\over{\ges}}\right)}}
\def\jac{J_{\g}}
\overfullrule=0pt
\def\ref{\hangindent=\parindent \hangafter=1 \noindent} 


\centerline {\bigbf Computing Multidimensional Residues}
\vskip .7cm

\centerline{{\bf  Eduardo Cattani}}
\centerline{ Department of Mathematics}
\centerline{ University of Massachusetts}
\centerline{ Amherst, MA 01003, \ U.S.A. }

\vskip .4cm

\centerline{\bf Alicia Dickenstein}
\centerline{Dto. de Matem\'atica, F.C.E. y N. }
\centerline{Universidad de Buenos Aires}
\centerline{ Ciudad Universitaria - Pab. I }
\centerline{(1428) Buenos Aires, \ Argentina }

\vskip .4cm

\centerline{\bf Bernd Sturmfels}
\centerline{ Department of Mathematics}
\centerline{ Cornell University }
\centerline{ Ithaca, NY 14853, \  U.S.A. }

\vskip 1.4cm

{\bigbf Abstract }

Given $n$ polynomials in $n$ variables with a finite number of complex
roots, for any of their roots there is a local residue operator
assigning a complex number to any polynomial.  This is
an algebraic, but generally not rational,
function of the coefficients.
On the other hand, the global residue, which is defined
as the sum of the local residues over all roots,  has invariance properties
which guarantee its rational dependence on the coefficients
[9],[27]. In this paper we present symbolic algorithms for evaluating
that rational function.

Under the assumption that the
deformation to the initial forms is flat, for some choice of
weights on the variables,
we express the global residue as a single residue integral with respect to
the initial forms. When the input equations are a
Gr\"obner basis with respect to a term order, this leads
to an efficient series expansion algorithm for global residues,
and to a vanishing theorem with respect to the
corresponding cone in the Gr\"obner fan.

The global residue of a polynomial
is shown to equal the highest coefficient
of its (Gr\"obner basis) normal form, and, conversely,
the entire normal form is
expressed in terms of global residues.
This yields a new method for evaluating traces
over zero-dimensional complete intersections.
Applications to be discussed include the counting of real roots
(as in [4],[22]),  the computation of the degree of a polynomial map
(cf.~[12]), and the evaluation of multivariate symmetric
functions (cf.~[16],[21]). All results and algorithms are
illustrated for an explicit system in three variables.

\vfill\eject

{\bigbf 0.\ Basic Properties of Multidimensional Residues}
\medskip

Multidimensional residues play a fundamental role in
complex analysis and geometry.
Recent applications of residues in computational algebra
include explicit division formulae [5],[6],
the evaluation of symmetric functions [21],
the membership problem for polynomial ideals [9],[10],
the effective Nullstellensatz [13],
and numerical algorithms for solving polynomial systems [7].
In most of these articles the emphasis lies on
degree estimates and complexity results.
Our goal here is to develop practical tools
for computing global residues.
We thus refrain from using ``univariate projections''
or ``linear changes of coordinates''; instead we seek algorithms
involving Gr\"obner bases and sparsity-preserving series expansions.
While initially our discussion follows a path similar to  [25],[27],
it then proceeds to systematically develop
the interplay between residues and Gr\"obner bases.

This paper is organized as follows.
In \S 1 we consider $n$ polynomials in $n$ variables which form
a Gr\"obner basis (or H-basis) with respect to some choice of
positive weights. This hypothesis has natural geometric (1.3'),
algebraic (1.5) and analytic (1.7) interpretations.
In (1.17) we express the global residue as
a single residue integral with respect to the initial forms.
In \S 2 we specialize to the case of a Gr\"obner basis
in the usual sense, with respect to a term order.
In  (2.3) we express the residue as a coefficient of
a certain polynomial. This yields a
polyhedral vanishing theorem (2.5), and a bound on the degree
of the residue as a polynomial in the trailing coefficients (2.7).
The constructions of \S 1 and \S 2 lead to
algorithms, which will be presented in \S 3.
In \S 4 we relate  global residues to the
coefficients in the (Gr\"obner basis) normal form.
The global residue of a polynomial is shown to
equal the highest coefficient of its normal form (4.2).
This results in fast procedures
for computing residues and traces (4.8).
In \S 5 we present applications to symmetric functions, to computing
the degree of a polynomial map, and to
counting real roots. Finally, we study in \S 6 an explicit
system in three variables.

\vskip .1cm

In this section (\S 0) we review the complex-analytic definition and
basic properties of multidimensional residues. Details and proofs
can be found in [1],[15],[25]. For the algebraic counterpart to the
analytic theory see e.g.~[3],[19],[20],[23]. We shall discuss
the equivalence of the algebraic and the analytic approach
briefly at the end of \S 0.

Given $n$ holomorphic functions
$g_1,\ldots,g_n$ in an open set $U\subset \C^n$ with a single common zero $p$
in
$U$, one can associate to any holomorphic function $h\in \OO(U)$ the local
residue at $p$ of the meromorphic $n$-form
$$\omega = {{h(\x)\, d\x}\over
{g_1(\x)\ldots g_n(\x)}}\ ,\quad d\x= dx_1\wedge\ldots\wedge dx_n
\,.$$
This defines  the $\C$-linear operator
$$ \OO(U) \, \to \, \C \,,\quad
h \,\mapsto \,{\rm Res}_p\left({{h(\x)\, d\x}\over
{g_1(\x)\ldots g_n(\x)}}\right) :=
{{1}\over{(2\pi i)^n}} \int_{\Gamma_{\g}(\epsilon)} \omega\,,\leqno{(0.1)}$$
where $\Gamma_{\g}(\epsilon)$ is the real $n$-dimensional cycle
$$\Gamma_{\g}(\epsilon) = \{\,\x\in U \ :\ |g_i(\x)| = \epsilon_i\ ,
i=1,\ldots, n\,\}$$
with orientation defined by the $n$-form $d\,{\rm arg}(g_1)
\wedge\ldots\wedge d\,{\rm arg}(g_n)$,
and where $\epsilon = (\epsilon_1,\ldots,\epsilon_n)$ is any $n$-tuple of
sufficiently small, positive real numbers.

The Cauchy formula in $n$ complex variables provides the simplest example
of a residue operator:  If $g_i(\x) = (x_i - p_i)^{a_i+1}$, $a_i\in\N$,
then
$${\rm Res}_p\left({{h(\x)\, d\x}\over
{\prod_{i=1}^n (x_i - p_i)^{a_i+1}}}\right) \,\,=\,\,
{{1}\over{a_1!\ldots a_n!}} \left(
{{\partial^{a_1+\ldots +a_n}\,h}\over{\partial x_1^{a_1}\ldots
\partial x_n^{a_n}}}\right)(p)\,.\leqno {(0.2)}$$
Let $\jac := {\displaystyle \det \left({{\partial g_i}\over{\partial
x_j}}\right)}$ denote the Jacobian of $\g = (\cges)$.  Then,
$${\rm Res}_p\left(
{{h\,\jac\, d\x}\over{\ges}}\right) = \mu_{\g}(p)\,h(p)\,, \leqno{(0.3)}$$
where $\mu_{\g}(p)$ denotes the  intersection
multiplicity of $\g = (g_1,\ldots,g_n)$ at $p\, $ [15, p.~662ff.].
If $p$ is a simple root, hence $\jac(p)\not= 0$, then
$${\rm Res}_p\left(
{{h\, d\x}\over{\ges}}\right) = {{h(p)}\over{\jac(p)}}\,. \leqno{(0.4)}$$

Suppose now that $\cges\in \C[\x]$ are $n$-variate polynomials whose
zero set $Z(\g)$ is a non-empty finite subset of $\C^n$.  We can consider the
{\sl global residue operator} ([9],[15],[25]) which assigns to a
polynomial $h\in\C[\x]$ the complex number
$$\rg h = \rh := \sum _{p\in Z(\g)} {\rm Res}_p\left(
{{h\, d\x}\over{\ges}}\right)\,. \leqno{(0.5)}$$
The main functorial properties of the global residue are encapsulated in the
following two results whose proof may be found, for example, in [25, II.8.3-4].

\nnproclaim{(0.6) Transformation Law.}
Suppose that $\cfes\in\C[\x]$ have finitely many common roots and that we can
write $$f_i = \sum_{j=1}^n A_{ij} g_j \ ;\quad A_{ij}\in\C[\x]\ ,\quad
i=1,\ldots,n\,.$$
Then, for $h\in\C[\x]$,
$$\rh = {\rm Res}\left(
{{h\,\det(A_{ij})\, d\x}\over{\fes}}\right)\,.$$
\endproclaim
\bigskip
If $\cges\in\C[\x]$ are as above, the ideal $I$ generated by them is
zero-dimensional, and therefore
$V=\C[\x]/I$ is a finite dimensional $\C$-vector space.  Since  the
global residue $\rg h$ vanishes for $h\in I\,$ (see [15, p.~650]),
it defines a  $\C$-linear map:
$${\rm Res}_{\g}\colon V \to \C\ ,\quad h \mapsto  \rg h\ \,. $$
\nnproclaim{(0.7) Duality.} A polynomial $h\in \C[\x]$ lies in $I$ if and only
if  $\rg {fh} =0$
for all $f\in\C[\x]$.
\endproclaim
This duality law may be interpreted as follows:
Let $V^* = {\rm Hom}_{\C}(V,\C)$. Then the
pairing
$$ V\times V^* \to V^*\,,\quad (b,\phi)\to (b.\phi)\,,$$
where $(b.\phi)(b') = \phi(b\,b')$, makes $V^*$ into a $V$-module.  Statement
(0.7)  is equivalent to the assertion that the residue operator
${\rm Res}_{\bf g}$ is a generator
of $V^*$ as a $V$-module.
\medskip
We recall that $\rg h$ is a rational function,
with integral coefficients, in the
coefficients of $\cges$. It may be computed in simply exponential
time  with respect to $n$, the number of variables; see [9],[21],[27].
A general procedure, suggested in these articles,
is to find univariate polynomials $f_1(x_1),\ldots,
f_n(x_n)$ in the ideal generated by $\cges$ and to transform the global
residue applying (0.6).  The global residues with respect to $\cfes$ are then
computed as a sum of products of univariate global residues.  This
general procedure is often too slow for practical computations.
We seek more efficient algorithms for ``nice'' situations,
such as the case when $\cges$  are a \gb for a term order.
The theory for such nice situations
is to be developed in the next section.

In closing let us mention the relationship between this analytic definition of
global residue and the algebraic definitions: of course, they coincide.
Write for each $i= 1, \ldots, n$,
$$g_i(\y) - g_i(\x) \quad = \quad \sum_{j=1}^n g_{ij}(\y,\x) (y_j - x_j) \ .$$
and denote $\Delta := {\rm det} (g_{ij})$.
Let $U \subset \C^n$ be the union of relatively compact
open neighborhoods isolating each of
the points in $Z(\g)$, let $\epsilon$ be any $n$-tuple of small positive
real numbers, and define $\Pi_{\epsilon} := \{ \x \in U :
|g_i(\x)|~<~\epsilon_i \ , \forall i=1, \ldots , n\}$. For any holomorphic
function $h$ on $U$, one can deduce from (0.2) and (0.6) (cf.[25, \S 17])
the following integral representation known as Weil's formula [26]:
$$ h(\x) = {{1}\over{(2\pi i)^n}} \int_{\Gamma_{\g}(\epsilon)}
{{h(\y) \Delta(\y,\x) d\y} \over {\prod_{i=1}^n (g_i(\y) - g_i(\x))}}\,,
\quad \x \in \Pi_{\epsilon}. $$
As $|g_i(\x)| < |g_i(\y)| $ for any $i\,, \x \in \Pi_{\epsilon}$ and $\y \in
\Gamma_{\g}(\epsilon)$, the integrand may be expanded as a
multiple geometric series
$${{h(\y) \Delta(\y,\x)} \over
{\prod_{i=1}^n (g_i(\y) - g_i(\x))}}\,\, = \,\,
h(\y) \Delta(\y,\x)  \sum_{\alpha_i\geq 0}
\biggl(\prod_{i=1}^n{{ g_i(\x)^{\alpha_i}} \over
{ g_i(\y)^{\alpha_i + 1}}}\biggr)\,, $$
which converges uniformly on compact subsets of
$\Gamma_{\g}(\epsilon)~\times~\Pi_{\epsilon}$. As a result of term-by-term
integration, we deduce that (cf.[25,\S 5],[5],[6])
$$h(\x)\,\,= \,\,\rg {h ( \,\cdot\,)
\ \Delta(\,\cdot\,, \x)} +
\sum_{|\alpha| \geq 1 }
{\rm Res} \left( {{h(\y) \Delta(\y,\x) \, d\y}\over {\prod_{i=1}^n
g_i(\y)^{\alpha_i + 1}}}\right)\,  \g^{\alpha}(\x)\,, \quad
\forall \x \in \Pi_{\epsilon}. \leqno{(0.8)}$$
Note that the second summand on the right is in the ideal generated by
$\cges$ in ${\cal O}(\Pi_{\epsilon})$ and
that $ \rg {h (\,\cdot\,) \ \Delta(\,\cdot\,, \x)}$
depends polynomially on $\x$.  If, in addition, $h$ is
a polynomial, then the fact that $Z(\g)$ is contained in $\Pi_{\epsilon}$
plus the fact that local analytic membership is equivalent to local
algebraic membership ([24]), imply that
$$h(\x) = \rg {h(\,\cdot\,) \ \Delta(\,\cdot, \x)}
 \quad \hbox{ on the quotient ring }
\, V\,. \leqno{(0.9)}$$
In general, (0.8) does not provide a representation
of their difference as a polynomial linear combination of $\cges$.
Under the hypothesis (1.3) below,
it follows from the vanishing statement in (1.18) that
the series becomes a finite sum, giving an
effective division formula with remainder which involves computing only
finitely many global residues associated to powers of $\cges$.
In summary, the formula (0.9) is the algebraic version of the
integral representation. It proves
that the global residue we are considering coincides with the
``trace'' associated to $\Delta$ as in [19, Appendix F] and [13],
and with the Kronecker
symbol (i.e., the dualizing linear form associated to $\Delta$) as in [3].

\bigskip
\bigskip

{\bigbf 1.\ Gr\"obner Bases for a Weight Partial Order}
\medskip
Let $\K$ be any subfield of the complex numbers $\C$ and let
 $g_i(\x)$,  $i=1,\ldots,n$, be polynomials in  $S = \K[\x]$,$\ \x
= (x_1,\ldots, x_n)$ .  Let $\w =
(w_1,\ldots,w_n)\in \N^n$ be a positive {\sl weight vector}.
The {\sl weighted degree\/} of a monomial
$\x^{\a} = x_1^{a_1}\ldots x_n^{a_n}$ is
$$
\deg_{\w}(\x^{\a}) = \ip {\w}{\a} = \sum_{i=1}^n w_i\,a_i\,.$$
  We extend the notion of weighted degree to
arbitrary polynomials in $S$ in the usual manner.
Write each polynomial $g_i(\x)$ as
$$g_i(\x) = p_i(\x) + q_i(\x) \leqno {(1.1)}$$
where $p_i$ is $\w$-homogeneous, and
$$d_i=\dw {p_i} = \dw {g_i}\,;\quad \dw {q_i} < \dw {g_i}\,.\leqno
{(1.2)}$$

Throughout this section we make the following assumption:
$$ p_1(\x) =\cdots = p_n(\x) = 0 \quad\hbox{if and only if}\quad\x=0
\leqno{(1.3)}$$
In what follows we will interpret this condition geometrically
(1.3'), algebraically (1.5) and analytically (1.7).  Let
$$\tilde g_i(t;\x) = t^{d_i} g_i(t^{-w_1}x_1, \ldots,
t^{-w_n}x_n) \leqno{(1.4)}$$
This is a homogeneous polynomial of degree $d_i$ in
$(t;x_1,\ldots,x_n)$ relative to the weights $(1;w_1,\ldots,w_n)$.

Let $\P^n_{\w}$ denote the weighted projective space with homogeneous
coordinates $(t;x_1,\ldots,x_n)$ and weights $(1;w_1,\ldots,w_n)$.
The image of the hyperplane $\{t=1\}\subset\C^{n+1} \!\setminus\! \{0\}$
in $\P^n_{\w}$ is identified with $\C^n$. If
$$\tilde D_i = \{ (t;\x)\in \wp\ :\ \tilde g_i(t;\x) = 0\}$$
then (1.3) is equivalent to the geometric condition
$$\tilde D_1\cap\ldots\cap
\tilde D_n\subset \C^n \leqno {(1.3')}$$
\medskip
The algebraic meaning of (1.3) is best expressed using the
following notion of a
Gr\"obner basis: Given a polynomial $f\in S$, we denote by $\iw f$ its
form of highest weighted degree.  For any ideal $I\subset S$
we define the initial ideal $\iw I$ to be the ideal generated by
$\iw f$ where $f$ runs over $I$.  A finite subset ${\cal G}\subset I$
is said to be a {\it Gr\"obner basis} for $I$, relative to
the weight $\w$, provided:
$$\iw I \quad = \quad \langle \iw g\ :\ g\in {\cal G}\rangle\,.$$
We emphasize that $\iw I$ need not be a monomial ideal.
Some authors prefer to call ${\cal G}$ an {\sl H-basis}, a term which
goes back to Macaulay in the classical case $ \w = (1,1,\ldots,1) $.

\medskip
\nnproclaim {(1.5)\ Lemma. }  Suppose ${\cal G} = \{g_1,\ldots,g_n\}\subset
S$ satisfy (1.3). Then ${\cal G}$ is a Gr\"obner basis for the ideal it
generates.
Conversely, suppose ${\cal G} = \{g_1,\ldots,g_n\}\subset
S$ is a \gb, with respect to $\w$, for a zero-dimensional ideal $I$. Then
$\{g_1,\ldots,g_n\}$ satisfy (1.3).
 \endproclaim

\proof With the same notation as above we have $p_i = \iw {g_i}$ and
$q_i = g_i - \iw {g_i}$.  Since $p_1,\ldots,p_n$
define a complete intersection,
the Koszul complex on
these forms in exact. This implies that every syzygy
$\,\sum h_i \cdot e_i \, $ on $\,(p_1,\ldots,p_n)$
can be written as a linear combination of the basic
syzygies $\,p_j \cdot e_k \,- \,p_k \cdot e_j $.

Suppose that ${\cal G}$ is not a Gr\"obner basis.
Then, there exists a polynomial
$$ f \quad = \quad \sum h_i g_i \quad = \quad
\sum h_i p_i \, + \, \sum h_i q_i $$
whose initial form does not  lie in $\langle p_1,\ldots,p_n \rangle $.
Hence $\, \sum h_i p_i  \, = \,0 $. By the remark above,
we can write $\,\sum h_i\cdot e_i \, = \,
\sum_{j,k} b_{jk} \cdot (p_j \cdot e_k \,- \,p_k \cdot e_j )$,
and the leading term of
$$\, \sum h_i g_i \  = \
\sum h_i q_i \  = \
\sum_{j,k} \,b_{jk} \cdot (p_j q_k \,- \,p_k q_j ) $$
must lie in $\langle p_1,\ldots,p_n \rangle $.
This is  a contradiction, completing the proof of the first statement.
To prove the converse, it suffices to note  $\dim \,(I)
= \dim\, (\iw I)=0 \, $ (see e.g.~[17]).
$\ \diamond $
\medskip
{\bf  (1.6)\ Remarks:} \  (i) The first part of Lemma
(1.5) remains true for any set of polynomials ${\cal G} =
\{g_1,\ldots,g_m\}\subset S$  whose initial forms define a complete
intersection.
\smallskip
(ii)  The initial ideal $\iw I$ is a flat deformation of the given ideal $I$
\ (see e.g.~[11, Ch.~6]).
\smallskip
(iii) The results in this section can be extended to
fields other than the complex numbers using the
deformation techniques in [20].

\bigskip

For each $t\in \C$, consider the map
$\g_t\colon\C^n\to\C^n$ defined by
$ \g_t(\x) = (\tilde g_1(t;\x),\ldots,\tilde g_n(t;\x))$.
We have the following analytic interpretation of (1.3).  Recall that  a
map $F\colon\C^n\to\C^n$ is said to be ${\sl proper\/}$ if the
inverse image of any compact set is compact.

\medskip
\nnproclaim{(1.7) Lemma.}  The polynomials $g_1,\ldots,g_n$ satisfy
condition (1.3) if and only if the map $\g_t$ is proper for every $t\in \C$.
\endproclaim
\proof  At $t=0$ we have $\g_0 = (p_1,\ldots,p_n)$. Since the polynomials $p_i$
are weighted homogeneous, the inverse image $\g_0^{-1}(0)$ is compact if and
only if $\g_0^{-1}(0)=\{0\}$.  Thus, if $\g_0$ is proper, then condition (1.3)
is satisfied.

For the converse it is enough to show that the map $\g$ is proper, since (1.3)
is a condition on just the initial form of the polynomials.  Let
$\tilde{\g}\colon\C^{n+1}\setminus\{0\} \longrightarrow \C^{n+1}\setminus\{0\}$
be defined by $$\tilde {\g}(t,x_1,\ldots,x_n) = (t,\tilde
g_1(t;\x),\ldots,\tilde g_n(t;\x))$$ The fact that $\g(\x)$ satisfies  (1.3)
guarantees that $\tilde {\g}$ takes values in $\C^{n+1}\setminus\{0\}$.
Since
$$\tilde{\g}(\lambda\, t,\lambda^{w_1}\,x_1,\ldots,\lambda^{w_n}\,x_n) =
(\lambda\, t,\lambda^{d_1}\,\tilde g_1(t;\x),\ldots,\lambda^{d_n}\,
\tilde g_n(t;\x))\,,$$
$\tilde{\g}$ defines a map from
$\P^n_{\w}$  to weighted
projective space $\P^n_{\d}$ with weights $(1;d_1,\ldots,d_n)$.
We may now consider the embedding of $\C^n$ in  $\C^{n+1}\setminus\{0\}$ as
the hyperplane $\{t=1\}$.  Since $t$ has weight one in both $\P^n_{\w}$
 and $\P^n_{\d}$, the natural projection from $\C^{n+1}\setminus\{0\}$ to
$\P^n_{\w}$
or $\P^n_{\d}$ is a homeomorphism of the hyperplane $\{t=1\}$ to its image.
Thus, $\tilde {\g}$ is a continuous extension of $\g$ to appropriate
compactifications of $\C^n$. If $K\subset \C^n$ is compact then
$ {\g}^{-1}(K)$ is compact since it coincides with
$\tilde {\g}^{-1}(K)$.\ $\diamond$ \bigskip

{\bf (1.8) Examples: }  An important special case, to be investigated in
detail in \S 2, is that of $n$ polynomials $g_1,\ldots,g_n\in S$ with finitely
many common roots in $\C^n$ and such that they are a \gb with respect to some
term order $\prec$. We can choose a weight vector
$\w\in\N^n$ such that  ${\rm in}_{\prec} (g_i) = \iw {g_i}$.
Since the ideal generated by the initial monomials $\iw {g_i}$  is
zero-dimensional, we may assume without loss of generality that:
$$\iw {g_i} = \alpha_i\,x_i^{r_i+1} $$
for some $ \alpha_i \in \K \setminus \{0\}, $ and therefore they satisfy
(1.3).  Particular examples are:

(i) The term order $\prec$ is lexicographic order:
then $g_1,\ldots,g_n$ satisfy
$g_i(\x) = g_i(x_i,\ldots,x_n)$ and $g_i$ is monic in $x_i$.
This is the case studied in [9]; see also [21] for the subcase when
$g_i(\x) = g_i(x_i)$ are univariate polynomials.

(ii) The term order is defined as total degree with ties broken by
lexicographic order with $x_n>\ldots>x_1$.  This is the case studied in
[1, (21.3)] and [25, II.8.2]. A weighted variant, due to
A\u{\i}zenberg and Tsikh, is studied in [1, (21.5)].

\bigskip
In this section  we are interested in studying the global
residue $\rg h$, for a polynomial $h\in\C[\x]$, under the
hypothesis that $g_1(\x),\ldots, g_n(\x)$ satisfy (1.3).
This hypothesis makes it possible to reduce the computation of $\rg h$ to that
of residues involving only certain powers of the initial forms
$p_1(\x),\ldots,p_n(\x)$. This is the content of (1.20) below.
In fact, considering $t$ as a parameter, the idea that one
can recover the information from the deformation to the initial forms
is the core of the geometric interpretation of Gr\"obner bases
(see e.g.~[2]).

Since the map
$\g(\x) = (g_1(\x),\ldots, g_n(\x))$ is proper,  we can replace (0.5) by a
single integral
$$\res
\left({{h\,d\x}\over{g_1\ldots g_n}}\right) \,\,= \,\,{{1}\over{(2\pi i)^n}}\,
\int_{\Gamma(\r)} {{h\,d\x}\over{g_1\ldots g_n}}$$ for any $\r=(r_1,\ldots,r_n)
\in (\R_{>0})^n$, where $\Gamma(\r)$  is the compact, real
$n$-cycle $$\Gamma(\r) = \{\,\x\in\C^n\ : \ |g_i(\x)|=r_i\ ;\ i=1,\ldots
n\,\}\,.$$

Similarly, for each fixed $t\in \C$,
$$R_h(t) := \res \left({{h(\x)\,d\x}\over{\tilde g_1(t;\x)\ldots \tilde
g_n(t;\x)}}\right) =   {{1}\over{(2\pi i)^n}}\
\int_{\tilde \Gamma_t(\r)} \,{{h(\x)\,d\x}\over{\tilde g_1(t;\x)\ldots \tilde
g_n(t;\x)}}\,,\leqno(1.9)$$
where the global residue is taken
relative to the divisors $
\left\{\,\x\in\C^n\ :\ \tilde g_i(t;\x) = 0\, \right\}$,
$\,i=1,\ldots,n\,$ and $\tilde \Gamma_t(\r) =
\{\,\x\in\C^n\ : \ |\tilde g_i(t;\x)|=r_i\ ;\ i=1,\ldots n\,\}
$.

The family $\{ \tilde g_1(t,\x), \dots , \tilde g_n(t,\x) \}$ is a \gb
with respect to $\w'=(1,2w_1, \ldots , 2w_n)$ for the ideal $\tilde I$
generated by $\{ \tilde f , f \in I \}$ and ${\rm in}_{{\w'}}
(\tilde f) = {\rm in}_{\w}(f) , \ \forall f \in S. $
It then follows that the coordinates
$(t,\x)$ are in {\sl Noether position} [8] for $\tilde I$ and, consequently,
we may apply Theorem 3.3 in [9] to deduce that $R_h(t)$ is a polynomial in $t$.
We will reprove this and in fact obtain
the stronger result (1.18).

We set
$$\tilde G(t;\x) \,\,\, := \,\,\, \prod_{i=1}^n \tilde g_i(t;\x)
\,\,\, = \,\,\, \sum_{j=0}^{d_{\w}}
A_j(\x)\,t^j \leqno {(1.10)}$$
where $d_{\w} = d_1+\ldots+d_n$.  Inverting the polynomial  $\tilde G(t;\x)$
as a rational formal power series in $t$, we write
$$\tilde G^{-1}(t;\x)\  = \
\sum_{j\geq 0} B_j(\x)\,t^j\,. \leqno {(1.11)}$$
Given positive real numbers $k_1,\ldots,k_n$, let
$$T(\k): = \{\,\x\in\C^n\,:\,|p_i(\x)| = k_i\,\}\,.$$
\bigskip
\nnproclaim {(1.12) Lemma.}  Given $\delta>0$, there exist positive constants
$k_1,\ldots, k_n$ so that, for $|t|\leq \delta$,
$$R_h(t) =
{{1}\over{(2\pi i)^n}}\
 \sum_{m\geq 0}\left(\int_{T(\k)} \,h(\x)\,B_m(\x)\,d\x\right) \,t^m\,,$$
and this series is uniformly convergent for $|t|\leq \delta$.
\endproclaim
\proof
Suppose $g_1,\ldots,g_n$ satisfy (1.3) and let
$$\tilde g_i(t;\x) = p_i(\x) + t\,\hat q_i(t;\x).$$
Then, given $\delta > 0$, there exist positive constants $k_1,\ldots,k_n$ such
that, for
$\x\in T(\k)$,
$${{|p_i(\x)|} \over {2}} > |t\,\hat q_i(t;\x)|\leqno{(1.13)}$$
for all $i= 1,\ldots, n$ and $|t|\leq\delta$.
Indeed, because of the weighted homogeneity property of $p_i$, it suffices to
take $k_i = \lambda^{d_i}$ with $\lambda$ sufficiently large.

The estimate (1.13) allows us to apply
{\sl Rouch\'e's principle for residues}
[25, II.8.1] and replace, for $|t|\leq \delta$, the integration cycles
$\tilde \Gamma_t(\r)$ in (1.9) by the fixed cycle $T(\k)$.  Thus,
$$R_h(t) =   {{1}\over{(2\pi i)^n}}\
\int_{T(\k)} \,{{h(\x)\,d\x}\over{\tilde g_1(t;\x)\ldots \tilde
g_n(t;\x)}} \qquad \qquad \hbox{for all $|t|\leq \delta$.}
\leqno{(1.14)}$$
In view of (1.13), it follows that the series
$$\sum_{j\geq 0} B_j(\x)\,t^j \quad  = \quad
{ {1} \over
{\displaystyle{ \prod_{i=1}^n p_i(\x)} \left(1+{\displaystyle  {{\strut t\,\hat
q_i(t;\x)}\over{p_i(\x)}}}\right)} }$$
is uniformly convergent for $\x\in T(\k)$ and
$|t|\leq\delta$.  Since we can now integrate (1.14) term-by-term, the result
follows. \ $\diamond$

\bigskip
\nnproclaim {(1.15) Lemma.} Let $P(\x) = p_1(\x)\ldots p_n(\x)$.  Then
$P^{m+1}(\x)\,B_m(\x)$ is a weighted homogeneous polynomial of degree
$ \,\,m\,(d_{\w}-1) \,$ with respect to $\w $.
\endproclaim
\proof
Since $\tilde G(t;\x)$ is weighted homogeneous of degree $d_{\w}$ and $t$ has
weight $1$, the coefficients $A_j(\x)$ in  (1.10) are weighted
homogeneous of degree $d_{\w} - j$.  On the other hand, the series (1.11)
inverts (1.10). This implies the following recursion relations:
$$\sum_{j=0}^m A_j\,B_{m-j} = 0 \quad ,\quad m\geq 1 \leqno{(1.16)}$$
with initial conditions $\, A_0\,B_0 =1 \,\,$ and $\,\, A_0({\x}) = P(\x)$.
In particular,
$$P\,B_m \,\, =  \,\, - \sum_{j=1}^m A_j\,B_{m-j}\,,$$
and
$$P^{m+1}\,B_m  \,\,=  \,\,\sum_{j=1}^m A_j\,P^{j-1}\,(P^{m-j+1}\,B_{m-j})\,.$$
Assuming that (1.14) holds inductively with respect to $m$, we obtain
$$\dw {P^{m+1}\,B_m} \,\,=
\,\,(d_{\w}-j) + (j-1)\,d_{\w} + (m-j)\,(d_{\w}-1) \,\, = \,\,
 m\,(d_{\w}-1) \,.\eqno{\diamond}$$

\medskip
The following is the main result in this section.

\nnproclaim {(1.17) Theorem. }
For any monomial $ \x^{\a} = x_1^{a_1}\ldots x_n^{a_n}$, set
${s(\a) = \langle \w,\a \rangle - d_{\w} +
\sum_{i=1}^n  w_i}$.  Then
$$R_{\x^{\a}}(t)\,\, = \,\,\rg {\x^{\a}}\cdot
t^{s(\a)}\,,$$
and
$$
\rg {\x^{\a}}\,\,\,=\,\,\,
{{1}\over{(2\pi i)^n}}\ \int_{T(\k)} \,\x^{\a}\,B_{s(\a)}(\x)\,d\x\,.
$$
\endproclaim
\medskip
Before proving (1.17), we note the following weighted version of the
Euler-Jacobi Theorem [15, p.~671]. A more general toric version
was given by Khovanskii in [18].
 \nnproclaim {(1.18) Corollary.}  $R_h(t)$ is a polynomial in $t$ of
degree at most  $\dw h - d_{\w} + \sum w_i$, and
$$\res \left({{h \, d\x}\over{g_1\ldots
g_n}}\right) = 0
\quad\qquad \hbox{whenever \ \ $\dw h < d_{\w} - \sum w_i$}. $$
\endproclaim

We observe also that, under the current hypothesis (1.3), this corollary
implies that the terms in the series in (0.8) will vanish for
$\sum_{i=1}^n \,\alpha_i\,d_i > \dw h$.
\medskip
{\bf Proof of (1.17):}  We begin by noting that, as in the case with
unit weights [25, IV.20.1]:

\noindent {\bf (1.19) } {\sl
If $P$ and $Q$ are weighted homogeneous polynomials in
$\C[\x]$, and  $\dw P - \dw Q + \sum w_i \not= 0\,,$
 then the form
$\,\omega = \displaystyle{{{P(\x)\, d\x}\over{Q(\x)}}}\,$
is exact. }

Indeed, we find that
$\,\omega = (\dw P - \dw Q + \sum w_i)^{-1} \ d\sigma \,$,
where $$\sigma \quad =
\quad {{P(\x)}\over{Q(\x)}}\  \sum_{j=1}^n (-1)^{j-1} w_j\,x_j \,
dx_1\wedge\ldots\wedge\widehat{dx_j}\wedge\ldots\wedge dx_n\,.$$
The verification of this equality is a straightforward consequence of Euler's
formula for weighted homogeneous polynomials:
$$\sum_{j=1}^n w_j\, x_j\, {{\partial P}\over{\partial x_j}}\, =\,\dw P\, P
\,.$$
As in Lemma (1.12), we write
$$R_{\x^{\a}}(t) = \sum_{m\geq 0} \left(\int_{T(\k)}
\x^{\a}\,B_m(\x)\,d\x\,\right)\,t^m\,. $$
Since, by Lemma (1.15), $B_m(\x)$ is a quotient of weighted homogeneous
polynomials we can apply (1.19) to conclude that
$$\int_{T(\k)}
\x^{\a}\,B_m(\x)\,d\x\ = 0 $$
whenever
$$\dw {\x^{\a}} + \dw {B_m} + \sum w_i \quad \not= \quad 0\,. $$
This inequation is equivalent to $ m \not= s(\a)$.
Hence all integrals in (1.12) vanish, except for the one
with $ m = s(\a)$. This was precisely the claim of (1.17).
$\diamond$

The second assertion of (1.17) says  that we may write
$$
\res \left({{\x^{\a}\,d\x}\over{g_1\ldots g_n}}\right)\,=\,
{{1}\over{(2\pi i)^n}}\ \int_{T(\k)} \,
{{\x^{\a}\,(P^{s(\a)+1}(\x)\,B_{s(\a)}(\x))}\over{P^{s(\a)+1}(\x)}}\,d\x\,.
\leqno{(1.20)}$$
The numerator is a weighted homogeneous polynomial, by Lemma (1.15).
Therefore (1.20) is a residue with respect to the  $(s(\a)+1)$
power of the initial forms $p_1(\x),\ldots,p_n(\x)$.

We conclude this section by observing that as a direct consequence of (1.17)
and the Duality Theorem (0.7) we obtain \ (see [25, IV.20.1] for the
case of unit weights):

\nnproclaim {(1.21) Macaulay's Theorem.}
Let $p_1(\x),\ldots,p_n(\x)$ be weighted homogeneous polynomials whose only
common zero is the origin.  Then, any weighted homogeneous polynomial $h(\x)$
satisfying $$\dw h > d_{\w} - \sum_{i=1}^n  w_i $$
is in the ideal generated by $p_1(\x),\ldots,p_n(\x)$.
\endproclaim
\bigskip
\bigskip

{\bigbf 2.\ Gr\"obner Bases for a Term Order}
\medskip
In this section we  specialize to the case of $n$ polynomials
$g_1,\ldots,g_n\in S$ with finitely
many roots in $\C^n$, which are a \gb with respect to a
term order $\prec$.  As  in (1.8), we  choose a positive weight $\w\in\N^n$
 such that ${\rm in}_{\prec} (g_i) = \iw {g_i}$.
We may assume that
$$\iw {g_i} =  \ x_i^{r_i+1}\,,\quad i=1,\ldots,n\,. \leqno {(2.1)}$$
Let $\r = (r_1,\ldots,r_n)$ and $\rr = (r_1+1,\ldots,r_n+1)$. Then, with
notation as in \S 1, $d_{\w} = \ip {\w}{\rr}$, and
$$ P(\x) = x_1^{r_1+1}\ldots x_n^{r_n+1} = \x^{\rr}\,,$$
$$s(\a) = \ip {\w}{\a} - d_{\w} + \sum  w_i = \ip {\w}{\a - \r}\,,$$
and (1.15) implies the following homogeneity of the coefficients of
(1.11):
\medskip \nnproclaim {(2.2)\ } $B_m(\x)$ is a $\w$-homogeneous
Laurent polynomial of weighted degree $-(m+d_{\w})$.
\endproclaim
\medskip
Consequently, the integrand in (1.17) is a Laurent polynomial. Since
$$\int_{T(\k)} \,\x^{\b}\,d\x \,\,\,\,= \,\,\,0
\qquad \hbox{for $\,\,\b\not=(-1,\ldots,-1)$},$$
 we obtain the following result.
\medskip
\nnproclaim {(2.3) Theorem. }
The residue $\rg {\x^{\a}}$ equals
the $x_1^{-1}\!\ldots x_n^{-1}$-coefficient of  $\,\x^{\a}\,B_{s(\a)}(\x)\,$.
\endproclaim

\medskip
Theorem (2.3) is essentially a restatement of a formula for computing global
residues due to A\u{\i}zenberg and Tsikh [1, (21.5)], which,
in a simpler version pointed out to us by J.~Petean, says that $\rg {\x^{\a}}$
equals the $x_1^{-1}\ldots x_n^{-1}$-coefficient in the expression
$$
\sum_{|\alpha|\leq {\ip {\w}{\a -\r}}}\left( (-1)^{|\alpha|}\, \x^{\a-(\rr)}\,
\prod_{i=1}^n\,(q_i(\x)/x_i^{r_i+1})^{\alpha_i}\right)\,.$$
The introduction of the homogenizing parameter $t$ organizes the
computation of this Laurent series and the search for
the desired coefficient, as evidenced in Algorithm (3.1) below. Keeping
track of the homogeneity properties of the coefficients $B_j(\x)$, also
allows us to get more precise information about the global residues,
such as Theorems (2.5) and (2.7).

 {\bf (2.4) Remark.}  Note that (2.2) implies that
$\x^{\a}\,B_m(\x)$ may contain a term of the form  $\alpha\,x_1^{-1}\ldots
x_n^{-1}$ only if $$\ip {\w}{\a} - (m+d_{\w}) = - \sum w_i$$
that is, only if
$$ m = \ip {\w}{\a-\r}. $$
Combined with (1.12), this gives  a simpler proof of the first statement in
(1.17) in the case when
 $g_1,\ldots,g_n$ are a \gb with respect to some
term order.

A similar argument combined with the recursive relations
(1.16) makes it possible to improve on the vanishing
statement in Corollary (1.18). Let ${\cal W}^*$ denote the
polyhedral cone in $\R^n$ which is positively spanned by
all lattice points of the form $ \rr -\b$, where
$\x^{\b}$ runs over all monomials appearing
in the expansion of $\,g_1 (\x) g_2 (\x) \ldots g_n (\x) $.

\medskip
\nnproclaim {(2.5) Theorem.} The residue
$\rg {\x^{\a}}$ vanishes if $\a - \r$ lies outside the  cone  ${\cal W}^*$.
\endproclaim
\proof
The condition on $\b$ in the definition of ${\cal W}^*$
is equivalent to saying that $\x^{\b}$ appears
in one of the coefficients $A_j(\x)$, $j=1,\ldots,d_\w$,
in the expansion (1.10) of $\tilde G(t;\x)$ as a
polynomial in $t$. The recursion relations (1.16) imply that
$B_m(\x)$ consists of terms $ k_{\alpha} \x^{\alpha}$
where the $n$-tuples $\alpha \in \Z^n$ are in the translated cone
$-((\rr) + {\cal W}^*)$.
Indeed, if $k_{\alpha} \x^{\alpha}$ is a term in $B_m(\x)$,
then (1.16) implies that for some $j=1,\ldots,m$, the Laurent polynomial
$A_j(\x)\,B_{m-j}(\x)$ contains a term of the form $c_{\alpha} \x^{\alpha +
\rr}$
and therefore,
$$\alpha + \rr \quad =\quad \b + \beta , $$
where $\x^{\b}$ is a
monomial  in  $A_j(\x)$ and $\x^{\beta}$ is a
monomial  in  $B_{m-j}(\x)$.  Then
$$\alpha + \rr \quad = \quad (\b - (\rr))+(\beta + \rr)$$
and the assertion follows by induction on $m.$

Now, Theorem (2.3) implies that $\rg {\x^{\a}} = 0$
unless the coefficient $B_{s(\a)}(\x)$ contains a
term which is a non-zero multiple of
$x_1^{-(a_1+1)}\ldots x_n^{-(a_n+1)}$.  But this is possible only if
$-(\a + {\bf 1})\in -((\rr) + {\cal W}^*)$, or equivalently,
if $\,\a-\r\in {\cal W}^*$.\
$\ \ \diamond$

\bigskip

Theorem (2.3) may also be used to study the dependence of $\rg {\x^{\a}}$
on the coefficients of the polynomials $g_1,\ldots,g_n$.
We write
$$g_i(\x) = x_i^{r_i+1} - \sum_{j=1}^{\nu_i} \,c_{ij}\,\x^{\a_{ij}},
\leqno {(2.6)}$$
and let ${\cal W} \subset {\bf
R}^n$ denote the closed convex cone of all vectors $\w\in {\bf
R}^n$ such that
$$\ip {\w}{\a_{ij}} \leq \ip {\w}{(r_i+1)\,{\bf e}_i}\ ;\quad \hbox
{for all}\
i=1\ldots,n\,;\,\, j=1,\ldots,\nu_i.$$
This cone  is the polar dual of the cone ${\cal W}^*$
defined above. By assumption, ${\cal W}$ has non-empty interior.
Note that ${\cal W}$ is the cone in the {\sl Gr\"obner fan} of
$I$ corresponding to the given term order (see e.g.~[14, \S 3.1]).
We have the following result.

\nnproclaim {(2.7) Theorem.}  The residue
$\,\rg {\x^{\a}}\,$ is a polynomial function in the coefficients $c_{ij}$.
\break Its degree in the variable $c_{ij}$ is bounded above by
$${\min_{\w\in {\cal W}}}\  {{\ip {\w}{\a-\r}}\over {\ip
{\w}{(r_i+1)\,{\bf e}_i-\a_{ij}}}}  \leqno{(2.8)}$$
and the total degree in the variables $\c = (c_{ij})$ is bounded by
$${\min_{\w\in int({\cal W}) \cap \Z^n
}}\  \ip {\w}{\a-\r} . \leqno{(2.9)}$$
\endproclaim

\proof
Given $g_i(\x)$  as in (2.6),
its weighted homogenization with respect to a weight
$\w=(w_1,\ldots,w_n)\in\Z^n \cap {\cal W}$ , is given by
$$\tilde g_i(t;\x) = x_i^{r_i+1} - \sum_{j=1}^{\nu_i} \,c_{ij}\,t^{\ip
{\w}{(r_i+1)\,{\bf e}_i-\a_{ij}}}\, \x^{\a_{ij}}$$
where ${\bf e}_i$ denotes the  $i$-th unit vector.  Set
$\rho_{ij} =  (r_i+1)\,{\bf e}_i-\a_{ij}$.

According to (2.3),
$\rg {\x^{\a}}$ equals the $x_1^{-1}\ldots x_n^{-1}$-coefficient of
$\x^{\a}\,B_{s(\a)}(\x)\,$, which is equal to the coefficient of
$t^{\ip {\w}{\a-\r}}\,x_1^{-1}\ldots x_n^{-1}$
in the expansion of
$${{\x^{\a}}\over{\x^{\rr}}}\
\sum_{i_1,\ldots,i_n\geq 0}\left(\sum_{j=1}^{\nu_1} \,c_{1j}\,t^{\ip
{\w}{\rho_{1j}}}\, \x^{-\rho_{1j}}\right)^{\! i_1}\ldots \,\,
\left(\sum_{j=1}^{\nu_n} \,c_{nj}\,t^{\ip
{\w}{\rho_{nj}}}\, \x^{-\rho_{nj}}\right)^{\! i_n}\,.\leqno{(2.10)}$$
It is now clear that the residue  depends polynomially on the coefficients
$c_{ij}$ and that, for a given choice of $\w\in {\cal W}$, its degree in
$c_{ij}$ is bounded
above by $\ip {\w}{\a-\r}/\ip
{\w}{(r_i+1)\,{\bf e}_i-\a_{ij}}\,. $

To prove the bound in (2.9) it suffices to observe that if a monomial
$\c^{\bf k}$ appears in $\rg {\x^{\a}}$, then by (2.10)
$$
\sum_{i,j} k_{ij} \ip {\w}{\rho_{ij}} \; = \, \ip {\w}{\a-\r}.
$$
As $ \ip {\w}{\rho_{ij}} \geq 1$ for  $\w
\in int({\cal W}) $ integral and all
$i,j$,  the claim follows.
$\ \diamond$
\medskip
{\bf (2.11) Remarks :}
i) By the Cauchy-Schwarz inequality, the bound (2.8)
is minimized by vectors $\w\in {\cal W}$ which are as
far as possible from $\a-\r$
and as close as possible to $(r_i+1)\,{\bf e}_i-\a_{ij}$.

ii) If we are given a  system of polynomials $\cfes$
which satisfy (2.1), and supposing that only the constant coefficients
are perturbed, say
$$g_i(\x) \,\,\,= \,\,\, f_i(\x) - c_i \quad \quad i=1, \ldots , n \, ,$$
then the bound in (2.9) implies the following:
If a monomial $c_1^{k_1} \dots c_n^{k_n}$ appears in $\rg {\x^{\a}}$, then
$\,\ip {\w}{\k} \leq \ip {\w}{\a-\r} $, that is,
the weighted degree with respect to $\w$
of  $\rg {\x^{\a}}$ in $(c_1,\ldots,c_n)$ is
bounded by $\ip {\w}{\a-\r}$ for all integral vectors
$\w\in int({\cal W})$.

\bigskip
\bigskip

{\bigbf 3.\ Deformation Algorithms for Global Residues}
\medskip
Let $\{\cges\} \subset S$ be a Gr\"obner basis as in \S 2, and let
$\w\in\N^n$ such that (2.1)  holds. We have shown that for any polynomial
$$h \,=\,\sum_{\ip {\w}{\a}\leq d} c_{\a}\,\x^{\a} \quad \in \quad S\,,$$
the computation of the global residue
$$\rh = \sum_{\ip {\w}{\a}\leq d} c_{\a}\,{\rm Res}\left(
{{\x^{\a}\, d\x}\over{\ges}}\right)$$
can be reduced to a sum of residues (at the origin) with respect to the family
of monomials $\{\,x_1^{r_1+1}\!,\ldots,x_n^{r_n+1}\,\}$ and some of their
powers.  Theorem (2.3) gives the following algorithm for computing
all global residues  up to a given weighted degree.
Note that Algorithm (3.1) respects
possible sparsity of the input polynomials.

\nnproclaim {(3.1) Algorithm. }

{\bf Input:}  $\ \w\in\N^n$, $\cges\in S$ satisfying (2.1), and $ d\in \N$.
\hfill \break
{\bf Output:} The global residues $\rg {\x^{\a}}$,
for all $\a\in\N^n$ such that $\ip {\w}{\a}\leq d\,.$ \endproclaim

{\bf Step 1:}  Define the weighted homogenizations
$\tilde g_i(t;\x) =
t^{w_i(r_i+1)}\,g_i(t^{-w_1}x_1,\ldots,t^{-w_n}x_n)$.
\hfill \break
{\bf Step 2:} If $\ip {\w} {\a - \r} \leq 0$,
then $\rg {\x^{\a}} = 0$.
\hfill \break
{\bf Step 3:}  Set $d'= d - \ip {\w}{\r}$. Compute the Taylor polynomial
$\sum_{j=0}^{d'} B_j(\x)\,t^j$ of degree $d'$
 for $(\prod \tilde g_i(t;\x))^{-1}$.
\hfill \break
{\bf Step 4:} For each $\a$ such that $\ip {\w}{\r} \leq
\ip {\w}{\a}\leq d$, find the
coefficient of $\x^{-(\a+{\bf 1})}$ in the Laurent polynomial $B_{\ip {\w}{\a -
\r}}(\x)$. It equals $\rg {\x^{\a}}$.
\medskip
The verification that ${\cal G}=\{\cges\}$ is a \gb and, if so, the choice of a
compatible weight $\w\in\N^n$ may be accomplished in at most
$O(n^{n+2}\,m^{2n-1}\,d^{(2n+1)n})$ arithmetic operations, where $m$ is a bound
for the number of monomials in each $g_i$ and $d$ is a bound for their degrees.
This was shown in [14, \S 3.2].

Naturally,  Algorithm (3.1) will be most efficient
when it is known {\sl a priori} that  ${\cal G}$ is a
\gb and a ``small''  compatible weight is given.
However, even if the given equations ${\cal G}$ are not a Gr\"obner basis
for any $\w$, then Algorithm (3.1) still serves as a useful subroutine.
To illustrate this, we describe a general procedure for computing the
global residue associated to {\sl any} complete intersection zero
dimensional ideal:

\nnproclaim {(3.2) Algorithm. }

{\bf Input:}  Polynomials $\cges$ in $S$  whose ideal $I$ is zero-dimensional.
\hfill \break
{\bf Output:} The global residue $\rg {h}$, for any specified
polynomial $h \in S$.
\endproclaim

{\bf Step 1:} Choose a ``good'' term order ``$\prec$''.
\hfill \break
{\bf Step 2:} Starting with $\{\cges\}$, run the Buchberger algorithm
towards a Gr\"obner basis, until the current
basis of $I$ contains polynomials  $\, f_1,\ldots,f_n \,$ with
$\,\ini_{\prec} (f_i) = x_i^{r_i +1}$.
\hfill \break
{\bf Step 3:} By keeping track of coefficients during Step 2, we
obtain an $n \times n$-matrix
$A= (A_{ij})$ of polynomials such that
$f_i =\sum_{i=1}^n A_{ij} \, g_j$, for $ \, i=1,\ldots,n.$
\hfill \break
{\bf Step 4:} Compute the desired residue via the
following formula:
 $$\rh \quad = \quad \res  \left(
{h \, \det(A) \, d{\x}}\over {\fes}  \right). \leqno (3.3) $$

\bigskip
{\bf (3.4) Remarks. } (i) In Step 1 we may take ``$\prec$'' to be
{\sl optimal} in the precise sense of [14, \S 3.3].
Such a choice is possible at almost no extra computational cost if
the {\it Gr\"obner basis detection} procedure of [14, \S 3.2]
had been run beforehand to test applicability of (3.1).

(ii) The termination and correctness of Step 2
follows from $\,{\rm dim} (\ini_{\prec} (I))
\, = \, {\rm dim}(I) =0$.

(iii) The correctness of (3.3) is just
the Tranformation Law (0.6). In order the evaluate the
right hand side of (3.3), we may use either Algorithm (3.1),
in case many residues are desired, or the formula to be presented
in (4.2) below, in case only one residue is desired.

(iv) \ As shown in [9, Theorem 3.3], it is possible to find
polynomials $\cfes$ and $A_{ij}$, $1\leq i,j \leq n$ with degrees bounded by
$\,nd^{2n}+d^n+d$, where $d=\max(3,\max\{\deg(g_i)\})$.

(v) If the polynomials $\cges$ satisfy (1.3) relative to some weight $\w$,
but not necessarily (2.1), then the weighted version of the
following argument due to Tsikh [25, II.8.3] describes how to find
polynomials $\cfes$ in the ideal $I = \langle \cges \rangle$ such that
$$\iw {f_i} \quad = \quad x_i^{\rho}\quad\hbox{for some } \rho.\leqno{(3.5)}$$
Indeed, Macaulay's Theorem (1.21) implies that, if $\, \rho > d_{\w} -
(w_1+\ldots+w_n)$, then
$$ x_i^{\rho}\,\, \in \,\, \iw I \,=\, \langle p_1,\ldots,p_n\rangle .$$
This degree bound allows the use of linear algebra (over $\K$) to determine
polynomials $A_{ij}$, $\,1\leq i,j \leq n$, such that
$x_i^{\rho}= \sum_{i=1}^n A_{ij} \, p_j$.
Moreover, $\deg_{\w}(A_{ij}) = \rho\,w_i - \deg_{\w}(g_j)$.  Let
$$f_i \,\, := \,\, \sum_{i=1}^n A_{ij} \, g_j
\,\, = \,\, x_i^{\rho} + \sum_{i=1}^n A_{ij} \, q_j$$
Since  $\deg_{\w}\left(\sum\, A_{ij} \, q_j\right) < \rho\,w_i$,
(3.5) holds, and (3.1) or (4.2) are applicable.
Naturally, we may use the Buchberger Algorithm to organize
the computation of the $A_{ij}$, which leads to a version of
Algorithm (3.2) which is applied to the initial forms
$ p_1,\ldots,p_n $ rather than the entire equations
$ g_1,\ldots,g_n $. In summary, there is plenty of room
for experimentation~!

\bigskip
The fact that the coefficients $B_j(\x)$  are
weighted-homogeneous implies that we can fix the value of one of the variables,
say $x_n=1$, and obtain a non-homogeneous version of Algorithm (3.1):  Set
$x_n=1$ everywhere, and, for each $\,\a\,$ with $\ip {\w}{\a}\leq d$, find the
coefficient of
${\displaystyle {{x_1^{-(a_1 + 1)}\ldots x_{n-1}^{-(a_{n-1} + 1)}}}}$
in $B_{\ip {\w}{\a-\r}}(x_1,\ldots,x_{n-1},1) $.  We note that other terms of
the form
$k\,{\displaystyle {{x_1^{-(a_1 + 1)}\ldots x_{n-1}^{-(a_{n-1} + 1)}}}}$
may appear in different coefficients $B_j(x_1, \ldots,x_{n-1},1) $, but the
homogenizing parameter $t$ keeps track of the only one contributing to the
residue.

In the classical case  $\w =(1, \ldots ,1)$, we can apply an argument of
Yuzhakov [27] to give the following geometric interpretation of (3.1).
As in \S 1, we imbed $\C^n$ in $\P^n$; the $n$-form
${\displaystyle{{\x^{\a}\,d\x}\over{g_1(\x)\ldots g_n(\x)}}}$
may be extended to a meromorphic form $\Phi$ in $\P^n$, which has $\{t=0\}$ as
a polar divisor of order $\ip {\w}{\a-\r}+1$
if and only if $\ip {\w}{\a-\r} \geq 0$.
The global residue in $\C^n$ may now be expressed as a single residue at
a point at $\infty$: if $P = (0, \dots,0,1)$, then
$${\rm Res}_{\bf g} ({\x^{\a}}) \quad = \quad (-1)^n \,
 {\rm Res}_P \,\Phi \quad = $$
$$ {{1}\over{\ip {\w}{\a \!-\!\r}!}} {\rm Res}_{0} \!\!\left(
x_1^{a_1} \!\! \ldots
x_{n-1}^{a_{n-1}}\,.\,{{{\partial^{\ip {\w}{\a-\r}}}\over
{\partial t^{\ip {\w}{\a-\r}}}}} \!\! \left( {{1} \over
{{\prod_{i=1}^{n} \! \tilde g_i(x_1,\ldots,x_{n-1},1;t)}
}}\right)\! (0)\ dx_1\wedge\ldots\wedge dx_{n-1}\! \right)$$
which, in turn, may be seen to equal the coefficient of
${\displaystyle {{x_1^{-(a_1 + 1)}\ldots x_{n-1}^{-(a_{n-1} + 1)}}}}$
in \break $B_{\ip {\w}{\a-\r}}(x_1,\ldots,x_{n-1},1) $. With suitable
modifications, the same interpretation holds for arbitrary weights.

\bigskip
\bigskip

{\bigbf 4.\ Normal Forms}
\medskip

In this section we give a formula expressing the coefficients of the
(Gr\"obner basis) normal form in terms of global residues.
In particular, the global residue of a polynomial
equals the highest coefficient of its normal form (4.2). This leads to a fast
algorithm for computing residues as well as traces over a zero-dimensional
complete intersection.

Suppose $g_1,\ldots,g_n\in \K[\x]$ satisfy
(2.1) with $\iw {g_i} = x_i^{r_i+1}$. Then
$V = \K[\x]/I$ is an Artinian ring of $\K$-dimension
$(r_1+1)(r_2+1) \cdots (r_n+1)$. Abbreviating
$$\II \quad := \quad \{\,\ii=(i_1,\ldots,i_n)\in\Z^n\ :\ 0\leq i_j \leq
r_j\,,\,\, j=1,\ldots,n\,\},$$
the set of monomials  $\{\,\x^{\ii}\ :\ \ii\in \II\,\}$
is a $\K$-vectorspace basis of $V$.
Every polynomial $h\in \K[\x]$
has a unique  {\it normal form}
$$\NF( h) \quad = \quad \sum_{\ii\in\II} \,c_{\ii}(h)\, \x^{\ii} .
\leqno{(4.1)}$$
The scalars $\,c_{\ii}(h) \in  \K \,$ are uniquely defined by the
property that $h \equiv \NF( h)\ ({\rm mod\ } I)$.  They are computed
using the division algorithm modulo the Gr\"obner basis
$\{g_1,\ldots,g_n \}$.

\medskip

\nnproclaim {(4.2) Lemma.}  With the notation as above, every
polynomial $h\in \K[\x]$ satisfies
$$\rh \quad = \quad c_{r_1,\ldots,r_n}( h) $$
\endproclaim

\proof
By linearity of the residue operator,
$\, \rg h = \sum_{\ii\in\II}\,c_{\ii}
\, \rg {\x^{\ii}}$.  However, for $\ii\in\II$, $\ii\not= \r$, we have
$\ip {\w}{\ii-\r}<0$ and, consequently, $\rg {\x^{\ii}}=0$ by Corollary (1.18).
On the other hand, it follows from Theorem (1.17) that
$${\rm Res}\left(
{{\x^{\r}\, d\x}\over{\ges}}\right) \quad = \quad {\rm Res}\left(
{{\x^{\r}\, d\x}\over{x_1^{r_1+1}\ldots x_n^{r_n+1}}}\right)
\quad = \quad 1. $$
This proves Lemma (4.2).\ $\diamond$
\bigskip
We now show that all coefficients of the normal form may be computed using
residues: \medskip
\nnproclaim {(4.3) Theorem.}  Fix an order on the index set $\II$, and let $M$
be the symmetric $|\II|\times|\II|$-matrix defined by
$$M_{\ii\jj} \quad := \quad {\rm Res}\left(
{{\x^{\ii}\,\x^{\jj}\, d\x}\over{\ges}}\right)\,,\quad \ii,\jj\in\II\,.$$
Then $M$ is invertible and for any $h\in\K[\x]$,
$$\bigl(c_{\ii}(h)\bigr)_{\ii\in\II} \quad = \quad M^{-1}\cdot
\Bigl({\rm Res}\bigl(
{{h\,\x^{\jj}\, d\x}\over{\ges}}\bigr)\Bigr)_{\jj\in\II}$$
\endproclaim
\proof The Duality Theorem (0.7) implies that the symmetric bilinear form
$$ V\times V \,\to \, \K \,, \qquad
( h_1,  h_2) \, \mapsto \, {\rm Res}\left(
{{h_1\,h_2\, d\x}\over{\ges}}\right) \leqno{(4.4)} $$
is non-degenerate.  The symmetric matrix $M$ represents (4.4)
relative to the basis  $\{\,\x^{\ii}\, :\ \ii\in\II\,\}$.
Therefore $M$ is non-singular.  The second claim
follows from the fact that
$${\rm Res}\left(
{{\x^{\ii}\,h\, d\x}\over{\ges}}\right)
\,\, = \,\,\sum_{\jj\in\II} c_{\jj}(h)\, {\rm Res}\left(
{{\x^{\ii}\,\x^{\jj}\, d\x}\over{\ges}}\right) \,\,= \,\,
\sum_{\jj\in\II} M_{\ii\jj}\,c_{\jj}(h)\,.\ \diamond$$
\medskip
{\bf (4.5) Remark.} As observed in [3], if we choose the
lexicographical order in $\cal I \/ $,
then the matrix $M$ has the triangular form
$$
  M \quad =\quad \pmatrix{ &  &  &  & 1  \cr
               & { 0} &  & {\cdot} &   \cr
               &  & { \cdot} &  &  \cr
		     & {\cdot} &  & { *} &  \cr
		     1 &  &  &  &  \cr }
$$
Consequently, $\det M = \pm 1$ and it is easy to compute the
inverse of $M$.

\medskip
We have seen in \S 0 that the Duality Law may also be interpreted as saying
that the
global residue operator ${\rm Res}_{\bf g}$  is a generator of $V^*$ as a
$V$-module.  There is another element in this $V$-module which is of special
interest in computational algebra, namely,
the morphism ${\rm tr}\in V^*$ which assigns to a polynomial
$h$ the trace of the endomorphism of $V$ given by multiplication by $h$.
The {\it trace} can be computed by normal form reduction as follows:
$$
{\rm tr}(h) \quad = \quad
 \sum_{\ii\in\II} \,c_{\ii}\bigl( h\cdot \x^{\ii}\bigr) .\leqno{(4.6)}
$$
On the other hand, it is known (see e.g.~[23, Satz (4.2)]) that
the trace may be expressed in terms of the global residue:
$${\rm tr}(h) \quad = \quad \sum_{p\in Z(\g)} \mu_{\g}(p) h(p)
\quad  = \quad {\rm Res}  \left(
{{h\,J_{\g}\, d\x}\over{\ges}}\right)\,.\leqno{(4.7)}$$
This expression, together with  (4.2), gives the following
formula for computing the trace:

\nnproclaim {(4.8) Algorithm. } \ \
Compute the trace of an element $h \in V$ by
$\,\,{\rm tr}( h) \, = \, c_{r_1\ldots
r_n}(h \cdot J_{\g})$.
\endproclaim

Thus, to find the trace of an element of $V$ over $\K$,
it suffices to run a single normal form reduction. We found
Algorithm (4.8) to be quite efficient in practice.
Additional speed can be gained by simple tricks, such as
replacing $J_{\g}$ by  $\NF(J_{\g}) $ in Algorithm (4.8),
and by storing previously computed normal forms of monomials.

\medskip

 From Theorem 2.7 we can derive
bounds on the degree of the trace for
$\cges$ as in (2.6) which satisfy (2.1). Note that the
corresponding cone $\cal W$ in the Gr\"obner fan
has nonempty interior. For each value of the parameters $c_{ij}$,
let $Z_\c$ denote the zero set of $\cges$.

\nnproclaim {(4.9) Theorem. }   For any $h \in K[\x]$, the parametric trace
$${\rm tr}(h) (\c) \quad = \quad \sum_{p \in Z_\c} \mu_{\g}(p) \, h(p)
$$
is a polynomial function of $\c = (c_{ij})$ with degree bounded above by
${\displaystyle { \min_{\w\in int({\cal W})
\cap \Z^n }} \dw {h}}$.
\endproclaim
{\bf Proof :} Let $\jac (\x)$ be the Jacobian of $\cges$ with respect to
the variables $x_1, \ldots,x_n$. Given any polynomial $h$ and
\hbox{$\w\in int({\cal W})
\cap \Z^n$}, we know by the second bound in Theorem 2.7 that the degree
of ${\rm tr}(h)(\c) = \rg {h \cdot \jac}$ in the $\bf c$ variables
is bounded by \hbox{$\dw {h \cdot \jac} - \ip{\w}{\r}$}. As $\dw {\jac}  =
\ip {\w}{\r}$, the claim follows.$\ \diamond$
\medskip

\bigskip

{\bigbf 5.\ Real Roots, Degree, and Symmetric Polynomials}
\medskip

 We present three applications of the computation of global residues:
counting real roots using the trace form (following
Pedersen-Roy-Szpirglas [22], Becker-W\"ormann [4]),
computing the degree of a polynomial map (following Eisenbud-Levine [12]),
and evaluating elementary symmetric polynomials in a multivariate
setting (following classical work of Junker [16]).

\medskip
We assume as above that $g_1,\ldots,g_n\in \K[\x]$ satisfy
(2.1) with $\iw {g_i} = x_i^{r_i+1}$, and again let
$\II \  := \  \{\,\ii=(i_1,\ldots,i_n)\in\Z^n\ :\ 0\leq i_j \leq
r_j\,,\,\, j=1,\ldots,n\,\}.$
Suppose that $\K$ is a subfield of the real numbers $\R$ and let $T$ be
the symmetric $|\II|\times|\II|$-matrix $T$ defined over $\K$ by
$$ T_{\ii\jj} \quad := \quad {\rm tr}\left(
\x^{\ii}\,\x^{\jj} \,\right)\,,\quad \ii,\jj\in\II\,.
\leqno (5.1) $$
The following result is due to
Becker-W\"ormann and Pedersen-Roy-Szpirglas.

\nnproclaim {(5.2) Theorem. } {\rm ([4],[22])} \ \ \
The rank of $T$ equals the number of distinct
complex roots in $Z({\bf g})$.
The signature of $T$ equals the number of distinct
real roots in $Z({\bf g}) \cap \R^n$.
\endproclaim

Recall that the {\it signature} of $T$ equals
the number of positive eigenvalues minus the number of negative
eigenvalues \ (all eigenvalues of a real symmetric matrix are real).
As is pointed out in [22, Prop.~2.8],
the signature can be read off directly from (the number of
sign variations in) the characteristic polynomial of $T$.
A straightforward generalization of (5.2) states that, for any $h \in V$,
the signature of the matrix $\,T^h \, =\, \bigl(
{\rm tr} (\x^{\ii}\,\x^{\jj} \cdot h ) \bigr)_{\ii,\jj\in\II} \, $
equals the number of distinct real roots with $h > 0$
minus the number of distinct real roots with $h < 0$.
Algorithms (4.8) and (3.1) provide subroutines
for computing $T$ and hence for counting real zeros
of zero-dimensional complete intersections.

\medskip

Viewing now $(g_1,\dots , g_n)$ as a proper map $\g : \R^n \to \R^n$, its
{\it degree} is defined as
$${\rm deg}(\g) := \sum_{p \in \g^{-1}(q)} {\rm deg}_p (\g)$$
where $q$ is a regular value of $\g$ and ${\rm deg}_p (\g)$ is $\pm 1$
depending on whether $J_{\g}(p)$ is positive or negative. The degree is a
topological invariant of $\g$.

Let $M$ be the non-singular, symmetric matrix $M$ defined, as in (4.3), by
$$M_{\ii\jj} \quad := \quad {\rm Res}\left(
{{\x^{\ii}\,\x^{\jj}\, d\x}\over{\ges}}\right)\,,\quad \ii,\jj\in\II\,.$$
The following result is essentially contained in [12]; although the results
there are local, the passage to the global situation may be done as in [22].

\medskip
\nnproclaim {(5.3) Theorem. } The degree of $\g$ equals the signature
of $M$.
\endproclaim

We may apply Algorithm (3.1) or Lemma (4.2) to compute the matrix $M$
and, consequently, the degree of $\g$.

\medskip
For our third application we need to review the concept of
symmetric polynomials in a multivariate setting. This
theory is classical (see Junker [16], who refers to even
earlier work of MacMahon and Schl\"afli). It reappeared
in the recent computer algebra literature in [21].
Let $A = (\alpha_{ij})$ be an $N \times n$-matrix of indeterminates
over $\K$. The symmetric group $S_N$ acts on the polynomial ring
$\K [ \alpha_{ij}]$ by permuting rows of $A$. We are interested in the
invariant subring $\K [ \alpha_{ij}]^{S_N} $, whose elements are
called {\it symmetric polynomials}. It is known that
$\K [ \alpha_{ij}]^{S_N} $ is generated
by symmetric polynomials of total degree at most $N$, but,
in contrast to the familiar $n=1$ case, this
$\K$-algebra $\K [ \alpha_{ij}]^{S_N} $ is not free
for $n \geq 2$. An important set of generators are the
{\it elementary symmetric polynomials} $\, e_{\bf j}(A) $, which
are defined as the coefficients of the following
auxiliary polynomial in $ u_1, u_2, \ldots u_n $:
$$ \prod_{ i=1}^N \,
(1 + \alpha_{i1} u_1 + \alpha_{i2} u_2 + \cdots + \alpha_{in} u_n )
\quad = \,\, \sum_{j_1+\ldots+j_n \leq N} \!\!\!
e_{j_1,\ldots,j_n}(A) \cdot u_1^{j_1}
u_2^{j_2} \cdots u_n^{j_n} \leqno (5.4) $$
Another set of generators is given by the {\it power sums : }
$$ h_{\bf j} \quad := \quad
\sum_{i=1}^N \alpha_{i1}^{j_1}
\alpha_{i2}^{j_2} \cdots \alpha_{in}^{j_n},
\qquad \hbox{for} \,\,\,
{\bf j} = (j_1,j_2,\ldots,j_n) \,\in \, \N^n ,\,\,
j_1+\ldots+j_n \leq N \leqno (5.5) $$
Algorithms and formulas for writing the $e_{\bf j}$ in terms
of the $h_{\bf j}$ and conversely are studied in detail
by Junker [16]. One of his methods will be
presented in (5.8)-(5.9) below.

Returning to our zero-dimensional complete intersection,
let $N = {\rm dim}_\K (V)$ be the cardinality of the multiset
$Z({\bf g}) \subset \C^n \,$ (counting multiplicities).
We fix any bijection between the rows of $A = (\alpha_{ij})$
and $Z({\bf g})$. This defines
a natural $\K$-algebra homomorphism
$$\phi \,\, : \,\, \K [ \alpha_{ij}]^{S_N} \quad \rightarrow
\quad \K, \leqno (5.6) $$
where the indeterminate $\,\alpha_{ij} \,$ gets
mapped to the $j$-th coordinate of the
$i$-th point in $Z({\bf g})$.
Our objective is to evaluate the map $\phi$ using
only operations in $\K$. In particular, we
are interested in the problem of
evaluating the elementary symmetric polynomials $e_{\bf j}$
under $\phi$.

The punch line of our discussion is that it is easy
to evaluate the power sums via the trace:
$$ \phi (h_{\bf j} ) \quad = \quad {\rm tr} ({\bf x}^{\bf j}) .
\leqno (5.7) $$
Thus to compute (5.7) we use Algorithm (4.8).
We then proceed using the following method due to
Junker and MacMahon. Consider the image
of (5.4) under $\phi$,
$$ R( {\bf u} ) \quad = \quad
\prod_{p \in Z({\bf g})}
(1 + p_{1} u_1 + \cdots + p_{n} u_n )^{ \mu_{\bf g}(p) }
\quad = \,\,\,\, \sum_{{\bf j}}
\phi(e_{\bf j}) \cdot {\bf u}^{\bf j}  \leqno (5.8) $$
The polynomial $R({\bf u})$ is the {\it Chow form} of the zero-dimensional
scheme defined by $I$. In computer algebra it is also
known as the {\it U-resultant\/}.
 Following[16, pp.~233, Eq.~(4)],
the formal logarithm of (5.8) equals
$${\rm log} \bigl( R( {\bf u} ) \bigr) \quad = \quad
\sum_{d=1}^\infty  {(-1)^{d-1} \over d} \cdot
\sum_{|{\bf j}| = d} \,
{ d \choose {\bf j}}  \,\phi (h_{\bf j} ) \,{\bf u}^{\bf j}.
\leqno (5.9) $$
Here $|{\bf j}| = j_1 + \cdots + j_n $ and
$\, { d \choose {\bf j}} \, = \,
{ d ! \over  j_1 ! j_2 ! \cdots j_n ! } $. Using (5.7) and (4.8),
we can compute the formal power series (5.9) up to any desired
degree $d' $. We then formally exponentiate this truncated series
(using operations only in $\K$) to get
the Chow form (5.8) up to the same degree $d' $.
In order to determine (5.8) completely, which means to evaluate
all elementary symmetric polynomials, it suffices to expand (5.9)
up to degree $d' = N = {\rm dim}_{\bf K}(V)$.

\bigskip
\bigskip

{\bigbf 6. An Example}
\medskip
In this section we apply our results and algorithms
to the specific trivariate system:
$$ g_1 \, = \, \underline{ x_1^5 } + x_2^3 + x_3^2 - 1 , \quad
 g_2 \, = \, x_1^2 + \underline{ x_2^2 } + x_3 - 1 , \quad
 g_3 \, = \, x_1^6 + x_2^5 + \underline{x_3^3 } - 1 .
\leqno (6.1) $$
This example is taken from [14, Example 3.1.2],
where it served to illustrate the problem of {\sl Gr\"obner Basis Detection}.
Indeed, the polynomials $g_1,g_2,g_3 $ are a Gr\"obner basis,
namely, for the weight vector ${\bf w} = (3,4,7)$.
With respect to these weights, the initial monomials are
the pure powers underlined above. We see that, counting possible
multiplicities, the set $Z({\bf g})$ consists of $30$ points in $\C^3$.
Our basic problem is to evaluate the global residue
$$ {\rm Res}_{\bf g} (\x^{\bf a}) \quad = \quad
 {\rm Res}\left({{x_1^{a_1} x_2^{a_2} x_3^{a_3}\, d\x}\over
{g_1(\x) g_2(\x) g_3(\x)}}\right), \leqno (6.2) $$
for any non-negative integer vector $ {\bf a} = (a_1,a_2,a_3)$.

It is interesting to compare the relative efficiency of Algorithm (3.1)
and the \gb reduction method deduced from Lemma (4.2).
In step 1 of Algorithm (3.1) we compute the
weighted homogenizations
$$ \tilde g_1 \, = \,  x_1^5 + t^3 x_2^3 + t^7 x_3^2 - t^{15} , \quad
 \tilde g_2 \, = \,  x_2^2  + t x_3 + t^2 x_1^2 - t^8 , \quad
 \tilde g_3 \, = \,  x_3^3 + t x_2^5 + t^3 x_1^6 - t^{21} . $$
We then consider the expression
$$ { 1 \over
\tilde g_1 (t;\x)  \cdot \tilde g_2 (t;\x)  \cdot \tilde g_3 (t;\x) }
\leqno (6.3) $$
as a rational function in $t$, and
we compute its Taylor expansion $\, \sum_{j = 0}^d B_j(\x ) \, t^j \,$
up to some degree $d $ which exceeds $3 (a_1 - 4) + 4 (a_2 -1) + 7 (a_3 -2)$.
Here the coefficients $B_j(\x)$ are
${\bf w}$-homogeneous Laurent polynomials in
$x_1,x_2,x_3$; for instance,
$$ B_2(\x) \quad = \quad
  { 1 \over x_1^{10} x_2^{4}}
- { 1 \over x_1^3 x_2^4 x_3^3 }
+ { x_2 \over x_1^5 x_3^5 }
+ { x_2^3 \over x_1^{10} x_3^4 }
+ { x_3 \over x_1^{15} x_2^2 }
+ { x_2^8 \over x_1^5 x_3^9 }
+ {1 \over x_1^5 x_2^6 x_3 } . $$
Now set $\,j = 3 (a_1 - 4) + 4 (a_2-1) + 7 (a_3-2) $.
The desired residue (6.2) equals the coefficient of
$\, x_1^{-a_1-1} x_2^{-a_2-1} x_3^{-a_3-1} \, $
in the Laurent polynomial $\, B_j  \,\in\,
\Z[x_1,x_1^{-1},x_2,x_2^{-1},x_3,x_3^{-1}]$.

This Taylor expansion is a fairly space
consuming process since the polynomials $B_j (\x)$ grow quite
large. This is witnessed by the following table, which shows
the number of terms of $B_j(\x)$ for some values of $j $ between $2$ and $40$:
$$ \matrix{
j : & 2  & 5 & 10 & 15 & 20 & 25 & 30 & 35 & 40  \cr
\# :  & 7  &  41 & 216 & 569 & 1102 & 1803 & 2682 & 3744 & 4964  \cr} $$

On the other hand, the normal form method of Lemma (4.2) is
quite efficient for evaluating individual residues.
Let $I$ denote the ideal in $\Q[x_1,x_2,x_3]$
generated by (6.1). The quotient ring $V = \Q[x_1,x_2,x_3]/I $
is a $30$-dimensional $\Q$-vector space.
Every element $h \in V$ is uniquely represented by its normal
form $\, {\cal N}{\cal F}(h) \, $ modulo the reduction relations:
$$  x_1^5 \,\longrightarrow \,
 - x_2^3 - x_3^2 + 1 , \qquad
  x_2^2 \,\longrightarrow \, - x_1^2 - x_3 + 1 , \qquad
  x_3^3  \,\longrightarrow \,- x_1^6 - x_2^5  + 1 . \leqno (6.4)$$
By Lemma (4.2), the residue (6.2) is equal to the coefficient of
$ x_1^4 x_2 x_3^2$ in $\, {\cal N}{\cal F}(h) $.

For instance, for the Jacobian
$$ J(\x) \,\, = \,\, {\rm det} \,\bigl( {{\partial f_i} \over {\partial x_j}}
\bigr)
\,\, = \,\,
 18 x_1^5 x_2^2 - 24 x_1^5 x_2 x_3  - 25 x_1^4 x_2^4
 + 30 x_1^4 x_2 x_3^2 + 20 x_1 x_2^4 x_3 - 18 x_1 x_2^2 x_3^2 $$
it takes $10$ reductions modulo (6.4) to reach the normal form
$$ \eqalign{  {\cal NF}(J)   \,\, = & \quad
\underline{30 x_1^4 x_2 x_3^2}
 - 25 x_1^4 x_3^2 - 152 x_1^4 x_2 + 146 x_1^4 x_3
 - 251 x_1^3 x_2 x_3 + 83 x_1^3 x_3^2 + 16 x_1^4
\cr & + 229 x_1^3 x_2
 + 8 x_1^3 x_3 - 196 x_1^2 x_2 x_3 + 226 x_1^2 x_3^2 - 114 x_1 x_2 x_3^2
 - 73 x_1^3 + 240 x_1^2 x_2
\cr & + 34 x_1^2 x_3 + 254 x_1 x_2 x_3 - 62 x_1 x_3^2
 + 69 x_2 x_3^2 - 260 x_1^2 - 140 x_1 x_2 - 78 x_1 x_3
\cr & + 108 x_2 x_3
- 49 x_3^2 +  140 x_1 - 177 x_2 - 128 x_3 + 177 . \cr } $$
Indeed, we see that the coefficient of $x_1^4 x_2 x_3^2$ equals
$\, {\rm Res}_{\bf g}(J) \,= \, {\rm tr}(1) \, = \, {\rm dim}\,(V) \,= \, 30 $.
Here is a slightly more serious example: it takes
$ 62 $ reductions modulo (6.4), running less than two minutes
in MAPLE on a Sparc 2, in order to find the global residue
$$ {\rm Res}_{\bf g}( x_1^{15}  x_2^{15}  x_3^{15} ) \quad =
\quad   -258,756,707,658,424,020,014,953,731,203.   $$

We made the observation that the efficiency of the two methods is comparable
when computing all residues of the form ${\rm Res}_{\g} (\x^{\a})$ with
$\ip {\w}{\a} \leq d$ for some fixed $d$. This is the case, for example,
in the computation
of the matrix $M$ defined in \S 4. This is a symmetric, $30 \times 30$-matrix
whose computation using Algorithm (3.1) requires the knowledge of
$B_j (\x)$ for $j \leq 30$. Using MAPLE on a Sparc 2 these may
be obtained in $321$ seconds. It takes an additional $247$ seconds
to read off the desired $465$ coefficients. On the other hand, it takes $324$
seconds to build up the matrix $M$ using Lemma (4.2).
The signature of $M$ is zero, and hence so is the degree of the map
$\g : \R^3 \to \R^3$ by Theorem (5.3).

For further combinatorial analysis we may wish to
compute the two polyhedral cones in Section 2.
We first obtain the $3$-dimensional quadrangular cone
$$ \eqalign{
 {\cal W} \quad = \quad & \bigl\{\,(w_1,w_2,w_3) \in \R^3 \,\, : \,\,
 5 w_1 \geq 2  w_3 , \,  2 w_2 \geq  w_3 ,\,
 w_3 \geq 2 w_1 ,\,  3 w_3 \geq 5 w_2 \bigr\} \cr
\quad = \quad & {\rm pos} \, \bigl\{
(4,5,10),\,(1,1,2),\,(5,6,10),\,(2,3,5) \bigr\}. \cr}$$
The interior of ${\cal W}$ consists of all weight vectors
which select the underlined monomials in (6.1) to be initial.
The cone polar to ${\cal W}$ equals
$$ \eqalign{ {\cal W}^* \quad  = \quad &
{\rm pos} \, \bigl\{ ( 5,  0, -2 ), \,( 0,  2, -1), \,( -2,  0,  1 ),\,
( 0, -5,  3) \bigr\} \cr \quad = \quad
& \bigl\{ \,(a_1,a_2,a_3) \in R^3 \,\, : \,\,
 4 a_1 + 5 a_2 + 10 a_3 \geq 0, \,
  a_1 + a_2 + 2 a_3 \geq 0,\,\cr
 & \qquad \qquad \qquad\qquad \quad \,
  5 a_1 + 6 a_2 + 10 a_3 \geq 0,\,
  2 a_1 +  3 a_2 +  5 a_3 \geq 0\bigr\}. \cr} $$
By Theorem (2.5), the residue (6.2) vanishes whenever
$$(a_1,a_2,a_3) \quad \not\in \quad (4,1,2) + {\cal W}^*,
\qquad \hbox{or equivalently,} $$
$$ 4 a_1 + 5 a_2 + 10 a_3 < 41 \,\hbox{ or } \,
  a_1 + a_2 + 2 a_3 < 9 \,\hbox{ or } \,
  5 a_1 + 6 a_2 + 10 a_3 < 46 \,\hbox{ or } \,
  2 a_1 +  3 a_2 +  5 a_3 < 21 .$$
For instance, $(6,1,1)$ satisfies the first inequality
and therefore $\, {\rm Res}_{\bf g}(x_1^6 x_2 x_3) \, = \, 0 $.

\vskip .1cm

In Section 4 we have shown that the trace ${\rm tr}(h)$
of an element $h$ in $V = \Q[\x]/I$ can be computed easily
as the coefficient of $ x_1^4 x_2 x_3^2$ in
$\, {\cal N}{\cal F}(h \cdot J) $. Using this technique,
let us now analyze the zero set $Z({\bf g})$
with respect to multiple roots, real roots, etc...
We compute the symmetric, integer $30 \times 30$-matrix
representing the trace form $T$ as in (4.8). The largest
entry in $T$ appears in the lower right corner:
$$ {\rm tr}(  x_1^4 x_2 x_3^2 \cdot   x_1^4 x_2 x_3^2 ) \quad = \quad
{\rm Res}_{\bf g} \bigl( x_1^8 x_2^2 x_3^4 \cdot J({\bf x} )\bigr)
\quad = \quad  16,049,138,278. $$
The rank of the matrix $T$ equals $20$.
By Theorem (5.2), this is the number of
{\sl distinct} roots of ${\bf g}$.
The characteristic polynomial of $T$ has
$13$ positive real roots and $7$ negative real roots.
Therefore the signature of $T$ equals $6$, and this is the
number of distinct {\sl real } roots of ${\bf g}$.
It turns out that there are four rational roots,
and they account for all multiplicities:
the root $(1,0,0)$ has multiplicity 3,
the root $(0,1,0)$ has multiplicity 4,
the root $(0,0,1)$ has multiplicity 6,
while the root $(-1,1,-1)$ is simple.
The remaining $16$ roots, two real and $14$ imaginary,
are  all simple and they are conjugates over $\Q $.

\vskip .1cm

We finally come to the problem of computing the Chow form
$$ R(u_1,u_2,u_3)
\quad = \quad \prod_{(\alpha_1,\alpha_2,\alpha_3) \in Z({\bf g})}
(1 + \alpha_1 u_1 + \alpha_2 u_2 + \alpha_3 u_3 )^{ \mu_{\bf g}(\alpha)} .$$
Note that each of the three non-simple roots
appears with its multiplicity in this product.
The ${33 \choose 3} = 5,456$ rational coefficients of $ R(u_1,u_2,u_3)$
are the values of the elementary symmetric polynomials at the roots of
${\bf g} = (g_1,g_2,g_3)$. Following (5.9), (5.7) and using Algorithm (4.8),
we compute the following formal power series up to a chosen degree:
$$  \eqalign{ &
{\rm log} \bigl( R(u_1,u_2,u_3) \bigr) \,\,\, = \,\,\,
 {\rm tr}(x_1) u_1+{\rm tr}(x_2) u_2+{\rm tr}(x_3) u_3  \,
- {1 \over 2} \cdot \bigl( {\rm tr}(x_1^2) u_1^2 + 2 {\rm tr}(x_1 x_2) u_1 u_2
\cr
& + 2 {\rm tr}(x_1 x_3) u_1 u_3 + {\rm tr}(x_2^2) u_2^2 + 2 {\rm tr}(x_2 x_3)
u_2 u_3
+ {\rm tr}(x_3^2) u_3^2 \bigr)
\,\, + {1 \over 3} \cdot \bigl( {\rm tr}(x_1^3) u_1^3 + \ldots \cr
& = \quad
 5 u_2 - 5 u_3
+ 37 u_1 u_2 - 121 u_1 u_3
- {35 \over 2} u_2^2 + 106 u_2 u_3 - {485 \over 2} u_3^2
+ 17 u_1^3 - 74 u_1^2 u_2
\cr &
+ 177 u_1^2 u_3 - 172 u_1 u_2^2
+ 536 u_1 u_2 u_3 - 686 u_1 u_3^2
+ {185 \over 3} u_2^3 - 667 u_2^2 u_3
+ 1084 u_2 u_3^2
\, + \, \ldots \cr } $$
By formally exponentiating this series, we obtain the Chow form
$$ \eqalign{ & R(u_1,u_2,u_3) \,\,\, = \,\,\,
1 \, + \, 5 u_2 - 5 u_3 \, + \,
37 u_1 u_2 - 121 u_1 u_3 - 5 u_2^2 + 81 u_2 u_3 - 230 u_3^2
\, + \, 17 u_1^3 \cr &
 - 74 u_1^2 u_2 + 177 u_1^2 u_3 + 13 u_1 u_2^2
 - 254 u_1 u_2 u_3 - 81 u_1 u_3^2 - 5 u_2^3
- 112 u_2^2 u_3 - 596 u_2 u_3^2
\, + \, \ldots \cr } $$
and hence all elementary symmetric polynomials.
For instance, we see that $\, \sum \alpha_1 \beta_2 \gamma_3
\, = \, -254 $, where the sum is taken over all triples
of roots $(\alpha_1,\alpha_2,\alpha_3),\,
 (\beta_1, \beta_2, \beta_3)  , \,
(\gamma_1, \gamma_2,\gamma_3) \in Z({\bf g})$.

\bigskip
\bigskip

{\bf Acknowledgements:}  This project began during the 1992 NSF
Regional Geometry Institute at Amherst College.  We thank its
organizers for their hospitality and, most particularly, its Research
Director, David Cox.  We express
our gratitude to Adrian Paenza and Paul Pedersen
for their help and support, and to the Center for Applied
Mathematics of Cornell University for its hospitality
during the preparation of this paper.
E.~Cattani was partially supported by NSF Grant
DMS-9107323, A.~Dickenstein was partially supported
by UBACYT and CONICET, and B.~Sturmfels was
partially supported by NSF grants DMS-9201453,
DMS-9258547 (NYI) and a David and Lucile Packard Fellowship.
\bigskip
\bigskip
{\bf Note added in proof:}  During the MEGA 94 meeting we became aware
of the paper: [M. Kreuzer and E. Kunz: Traces in strict Frobenius algebras
and strict complete intersections. {\sl J. reine angew. Math.\/} {\bf 381}
(1987), 181-204].  Our assumption (1.3) is equivalent, by their Proposition
(4.2), to the statement that the $K$-algebra $V = K[\x]/I$ is a strict
complete intersection.  Consequently, Theorem (1.17), in the case
$s(\a)\leq 0$ (in particular the Euler-Jacobi Theorem (Corollary (1.18))
is contained in their Corollary (4.6) and Theorem (4.8) .
\bigskip
\bigskip
\bigskip
\bigskip

\noindent{\bigbf References}
\medskip
\ref [1]\ \ \
I. A. A\u{\i}zenberg and A. P. Yuzhakov: Integral representations and
residues in multidimensional complex analysis.  Translations of
Mathematical Monographs {\bf 58}. American Mathematical Society, 1983.

\ref [2]\ \ \ D.~Bayer and D.~Mumford: What can be computed in Algebraic
Geometry?,  In:
{\sl ``Computational Algebraic Geometry and Commutative Algebra''\/}
(eds. D.~Eisenbud, L.~Robbiano),
Proceedings Cortona 1991, Cambridge University Press, 1993, pp.~1--48.

\ref [3]\ \ \ E.~Becker, J.-P.~Cardinal, M.-F.~Roy and Z.~Szafraniec:
Multivariate Bezoutians, Kronecker Symbol and Eisenbud-Levine formula,
submitted to MEGA 94.

\ref [4]\ \ \ E.~Becker and T.~W\"ormann: On the Trace Formula for
Quadratic Forms,
Recent Advances in Real Algebraic Geometry and Quadratic Forms,
Proceedings of the RAGSQUAD Year, Berkeley 1990-1991,
W.B. Jacob, T.-Y. Lam, R.O. Robson (editors),
{\sl Contemporary Mathematics}, 155, pp. 271--291.

\ref [5]\ \ \ C. Berenstein and A. Yger:
Une formule de Jacobi et ses cons\'equences. {\sl Ann. scient. Ec. Norm.
Sup.\/} 4e s\'erie, {\bf 24} (1991) 369--377.

\ref [6]\ \ \ C. Berenstein and A. Yger:
Effective Bezout identities in ${\Q}[z_1, \dots , z_n]$. {\sl Acta Math.\/}
{\bf 166} (1991) 69--120.

\ref [7]\ \ \ J.-P.~Cardinal: Dualit\'e et algorithmes it\'eratifs pour
la r\'esolution de syst\`emes polynomiaux, Th\`ese Univ. Rennes I, Janvier
1993.

\ref [8]\ \ \ A. Dickenstein, N. Fitchas, M. Giusti and C. Sessa:
The membership problem for unmixed polynomial ideals is solvable in single
exponential time. {\sl Discrete Applied Math.\/} {\bf 33} (1991) 73--94.

\ref [9]\ \ \ A. Dickenstein and C. Sessa: An effective residual criterion for
the membership problem in $\C[z_1, \ldots , z_n]$. {\sl Journal of Pure and
Applied Algebra\/} {\bf 74} (1991) 149--158.

\ref [10]\ \ \ A. Dickenstein and C. Sessa:  Duality methods for the membership
problem,  In:
{\sl ``Effective Methods in Algebraic Geometry''\/}
(eds. T.~Mora, C.~Traverso),
Proceedings MEGA-90, Progress in Math. 94, Birkh\"auser, 1991, pp.~89--103.

\ref [11]\ \ \ D. Eisenbud: Commutative Algebra with a View toward Algebraic
Geometry.  To Appear.

\ref [12]\ \ \ D. Eisenbud and H. Levine: An algebraic formula for the degree
of a ${\cal C}^{\infty}$ map germ. {\sl Annals of Mathematics\/} {\bf 106}
(1977) 19--44.

\ref [13]\ \ \ N.~Fitchas, M.~Giusti and F.~Smietanski: Sur la complexit\'e
du th\'eor\`eme des z\'eros, Preprint, 1993.

\ref [14]\ \ \ P.~Gritzmann and B.~Sturmfels: Minkowski addition of polytopes:
Computational complexity and applications to Gr\"obner bases. {\sl SIAM J.
Discr. Math.\/} {\bf 6} (1993) 246--269.

\ref [15]\ \ \ P. A.~Griffiths and J.~Harris:
Principles of Algebraic Geometry,
Wiley-Interscience, New York, 1978.

\ref [16]\ \ \ F.~Junker: \"Uber symmetrische Funktionen von
mehreren Reihen von Ver\"anderlichen, {\sl Mathematische Annalen},
{\bf 43} (1893) 225--270.

\ref [17]\ \ \ M.~Kalkbrener and B.~Sturmfels: Initial
complexes of prime ideals, {\sl Advances in Math.}, to appear.

\ref [18]\ \ \ A.~G.~Khovanskii: Newton's polyhedron and the Euler-Jacobi
formula.  {\sl Uspekhi Mat. Nauk\/} {\bf 33}, no. 6 (1978) 245--246; English
transl.~in {\sl Russian Math. Surveys} {\bf 33} (1978).

\ref [19]\ \ \ E.~Kunz: K\"ahler Differentials, Advanced Lectures in
Mathematics, Vieweg Verlag, 1986.

\ref [20]\ \ \ E.~Kunz and R.~Waldi: Deformations of zerodimensional
intersection schemes and residues,
{\sl Note di Matematica} {\bf 11} (1991) 247--259.

\ref [21]\ \ \ P.~Pedersen: Calculating multi-dimensional
symmetric functions using Jacobi's formula,
Proceedings AAECC 9, (eds. H.F.~Mattson, T.~Mora, T.R.N.~Rao),
Springer Lecture Notes in Computer Science,
{\bf 539}, 1991, pp. 304--317.

\ref [22]\ \ \ P.~Pedersen, M.-F.~Roy and A.~Szpirglas:
Counting real zeros in the multivariate case, In:
{\sl ``Computational Algebraic Geometry''\/}
(eds.~F.~Eyssette, A.~Galligo), Proceedings MEGA-92,
Progress in Math. 109, Birkh\"auser, 1993, pp.~203--223.

\ref [23]\ \ \ G.~Scheja and U.~Storch: \"Uber Spurfunktionen bei
vollst\"andigen Durchschnitten, {\sl J.~Reine u.~Angewandte Mathematik}
{\bf 278/9} (1975) 174--190.

\ref [24]\ \ \  J.-P.~Serre: G.A.G.A, {\sl Annales de l'Institut Fourier VI}
(1956),1--42.

\ref [25]\ \ \ A. K. Tsikh: Multidimensional residues and their applications.
Translations of Mathematical Monographs {\bf 103}. American
Mathematical Society, 1992.

\ref [26]\ \ \ A.~Weil: L'int\'egrale de Cauchy et les fonctions de
plusieurs variables. {\sl Math. Ann.} {\bf 111} (1935), 178--182.

\ref [27]\ \ \ A. P. Yuzhakov: On the computation of the complete
sum of residues relative to a polynomial mapping in $\C^n$.
{\sl Dokl. Akad. Nauk. SSSR\/} {\bf
275} (1984), 817--820;  English transl. in
 {\sl Soviet.~Math.~Dokl.\/} {\bf 29(2)} (1984) .

\bye